\documentclass[11pt]{amsart} 
\newcommand{\Hom}{\text{Hom}}
\newcommand{\Com}{\text{Com}}

\newcommand{\EEnd}{\mathcal End}
\newcommand{\Cup}{\text{Cup}}
\newcommand{\dev}{\text{dev}}
\newcommand\ad{\text{ad}}
\newcommand\bul{\bullet}
\newcommand{\1}{\text{id}}
\newcommand{\x}{\times}
\newcommand{\NN}{\mathbb{N}}
\newcommand{\I}{\hbox{I}}
\newcommand{\de}{\delta}
\newcommand{\la}{\lambda}
\newcommand{\EE}{\mathcal E}
\renewcommand{\=}{:=}
\renewcommand{\t}{\otimes}
\renewcommand{\o}{\circ}

\renewcommand\u{\smile}
\renewcommand\char{\text{char}} 
\begin{document}
\null
\vskip3true cm

\begin{center}
{\large\bf INVITATION TO COMPOSITION}\\

\vskip22pt
{\large Liivi Kluge}
\vskip5pt
Dept. of Mathematics,
Tallinn Technical University\\
Ehitajate tee 5, Tallinn 19086, Estonia\\
e-mail: liivi.kluge.ttu@mail.ee

\vskip15pt
{\large Eugen Paal}
\vskip5pt
Dept. of Mathematics,
Tallinn Technical University\\
Ehitajate tee 5, Tallinn 19086, Estonia\\
e-mail: epaal@edu.ttu.ee

\vskip15pt
{\large Jim Stasheff}
\vskip5pt
Dept. of Mathematics,
University of North Carolina\\
Chapel Hill, NC 27599-3250, U.~S.~A.\\
e-mail: jds@math.unc.edu
\end{center}

\begin{quote}
\item
\subsection*{Abstract}
In 1963 [Ann. of Math. {\bf 78}, 267-288], Gerstenhaber invented a 
\emph{comp(osition)} calculus in the Hochschild complex of an associative 
algebra. 
In this paper, the first steps of the Gerstenhaber theory are exposed
in an abstract (comp system) setting.
In particular, as in the Hochschild complex, a graded Lie algebra and a 
\emph{pre-coboundary} operator can be associated to every comp system.
A \emph{derivation deviation} of the pre-coboundary operator over the total
composition is calculated in two ways, (the long) one of which is essentially
new and can be seen as an example and elaboration of the
\emph{auxiliary variables method} proposed by Gerstenhaber in the 
early days of the comp calculus.

\smallskip
{\bf Classification.} 18D50~(MSC2000).

{\bf Key words.}
   Comp(osition), (pre-)operad, Gerstenhaber theory, cup, 
   graded Lie algebra, (pre-)coboundary, derivation deviation.
\end{quote}

\renewcommand{\thefootnote}{}
\footnote{\small 
  L.~Kluge and E.~Paal were supported in part by the ESF grants 1453 and 3654, 
  J.~Stasheff was supported in part by the NSF grant DMS-9504871.}

\baselineskip 1.5 \baselineskip
\section{Introduction and outline of the paper}

\noindent
In 1963, Murray Gerstenhaber \cite{Ger} invented a \emph{comp(osition)} calculus in 
the Hoch\-schild complex \cite{Hoch} of an associative algebra. 
The theory proposed in \cite{Ger} was announced \cite{Ger68,GGS92Am,CGS93}
to hold also in an abstract setting, i.~e. for abstract 
comp systems.
In this paper, the first steps of the Gerstenhaber theory 
(see also \cite{GerVor94,GetzJon94,VorGer,Akm97} for recent 
expositions and elaborations) are presented in an abstract comp system 
setting.
In particular, as in the Hochschild complex, a graded Lie algebra 
and a \emph{pre-coboundary} operator can be associated to
every comp system.
A \emph{derivation deviation} of the pre-coboundary operator over the total
composition is calculated in two ways.
The short way, however modified (simplified) in this paper, can be 
adapted from \cite{Ger} with some effort, 
but the long one (via auxiliary variables) 
is essentially new  and can be seen as an example and elaboration of the
\emph{auxiliary variables method} proposed by Gerstenhaber in the early days 
of the comp calculus.
We cover all the main aspects of \cite{Ger} from the
modern abstract point of view, except for the proof of Theorem~5 therein, 
because this needs quite a specific (auxiliary variables) technique, 
explained concisely in section \ref{aux}. 

As a variation of previous terms, we shall introduce notions of a 
\emph{composition algebra} and a $\u$-\emph{algebra}.
It turns out that \emph{right} translations in the composition algebra
are (right) derivations of the $\u$-algebra.
As may be anticipated from the Gerstenhaber theory of \emph{endomorphism}
comp systems, \emph{left} translations in 
the composition algebra are not derivations of the $\u$-algebra.
The corresponding \emph{derivation deviation} coincides up to sign with 
the derivation deviation of the pre-coboundary operator over 
Gerstenhaber's ternary braces (cf.~Theorem~\ref{third main} 
of our paper with Theorem~5 of \cite{Ger}).

\section{Pre-operad (composition system)}

\noindent
Let $K$ be a unital commutative associative ring, and let
$C^n$ ($n\in\NN$) be unital $K$-modules.
For \emph{homogeneous} $f\in C^n$, we refer to $n$ as the \emph{degree} 
of $f$ and write $(-1)^f\=(-1)^n$. Also, it is convenient
to use the shifted (\emph{desuspended}) degree $|f|\=n-1$.
Throughout this paper, we assume that $\t\=\t_K$.

\subsection{Definition {\rm (cf.~\cite{Ger,Ger68,GGS92Am,CGS93})}}
\label{pre-operad}

A linear (right) \emph{pre-operad} (\emph{composition system}) with 
coefficients in $K$ is a sequence $C\=\{C^n\}_{n\in\NN}$ of unital 
$K$-modules (an $\NN$-graded $K$-module), 
such that the following conditions hold.
\begin{enumerate}
\item 
  For $0\leq i\leq m-1$ there exist \emph{partial compositions} 
$$
  \o_i\in\Hom\,(C^m\t C^n,C^{m+n-1}), \qquad|\o_i|=0.
$$
\item
  For all $h\t f\t g\in C^h\t C^f\t C^g$, 
  the \emph{composition relations} hold,
$$
(h\o_i f)\o_j g=
   \begin{cases}
           (-1)^{|f||g|} (h\o_j g)\o_{i+|g|}f, 
                                     &\text{if $0\leq j\leq i-1$}\\ 
           h\o_i(f\o_{j-i}g), &\text{if $i\leq j\leq i+|f|$}.
\end{cases}$$
\item
   There exists a unit $\I\in C^1$ such that 
$$
\I\o_0 f=f=f\o_i \I,\qquad 0\leq i\leq |f|. 
$$
\end{enumerate}

\subsection{Remark} 

A pre-operad is also called a 
              \emph{comp(osition) algebra} 
           or \emph{asymmetric operad} 
           or \emph{non-symmetric operad} 
           or \emph{non-$\Sigma$ operad}. 
The concept of (\emph{symmetric}) \emph{operad} was formalized by 
Peter May \cite{may72} as a tool for the theory of iterated loop spaces.
Recent studies and applications can be found in \cite{Rene}.

Above we modified the Gerstenhaber \emph{comp algebra} defining relations 
\cite{CGS93,GGS92Am} by introducing the sign $(-1)^{|f||g|}$ in the 
defining relations of the pre-operad. 
The modification enables us to keep track of (control) sign changes more 
effectively.
One should also note that (up to sign) our $\o_i$ is Gerstenhaber's $\o_{i+1}$
from \cite{GGS92Am,CGS93}; we use the original (non-shifted) 
convention from \cite{Ger,Ger68}.

\subsection{Graded operads \cite{Markl,KiStaVo95}}
\label{graded operad}

The above definition makes sense also for \emph{internally} graded $C^n$ 
(i.~e. bigraded $C$), where $\o_i$ are of internal degree zero, $|\o_i|=0$
and the signs are adjusted as usual.
In this case, $|f|$ means the \emph{internal} degree of $f$, the convention
$|f|=f-1$ is abandoned and $\o_{i+|g|}\=\o_{i+g-1}$.
Such a pre-operad is called \emph{graded}.
The theory of \emph{differential graded operads} was recently developed in 
\cite{Markl,MarklShn94,MarklShn96,MarklShn96march}.

\subsection{Endomorphism pre-operad \cite{Ger,Ger68,GGS92Am}} 
\label{HG}

Let $A$ be a unital $K$-module and 
$\EE_A^n\={\EEnd}_A^n\=\Hom\,(A^{\t n},A)$.
Define the partial compositions for $f\t g\in\EE_A^f\t\EE_A^g$ as
$$
f\o_i g\=(-1)^{i|g|}f\o(\1_A^{\t i}\t g\t\1_A^{\t(|f|-i)}),
         \qquad 0\leq i\leq |f|.
$$
Then $\EE_A\=\{\EE_A^n\}_{n\in\NN}$ is a pre-operad 
(with the unit $\1_A\in\EE_A^1$) called the \emph{endomorphism pre-operad} 
of $A$.
A few examples (without the sign factor) can be found in
\cite{Ger68,GGS92Am} as well.
We use the original indexing of \cite{Ger,Ger68} for the
defining formulae. 

\subsection{Associahedra}

A geometrical  example of a pre-operad is provided by the Stasheff 
\emph{associahedra}, which was first constructed in \cite{Sta63}.
Quite a surprising realization of the associahedra as 
\emph{truncated simplices} was discovered and studied recently in 
\cite{ShnSte94,Sta97,Markl97}.
Unfortunately, the example is too sophisticated to expose in the 
present paper and we leave a reader alone with the above references.

\subsection{Representations}

It is widely accepted that abstract groups and algebras can be 
\emph{faithfully} represented by linear transformations.
In the theory of operads, \emph{endomorphism} operads play 
the same role with respect to \emph{abstract} ones.
Following this, one can expect that operads can be faithfully 
represented by endomorphism operads.
\emph{Representation} means that there exists an \emph{operad map}    
$\Psi\in\Hom\,(C,{\EEnd}_A)$ of degree \emph{zero}, such that
$$
\Psi(f\o_i g)=(\Psi f)\o_i(\Psi g), \qquad i=1,\cdots,|\Psi f|=|f|.
$$
The resulting triple $(C,A,\Psi)$ is called an \emph{algebra over the operad}
$C$ or $C$-\emph{algebra} in short.
In view of this, much of the Gerstenhaber theory for \emph{endomorphism} 
pre-operads is expected to hold for \emph{abstract} ones as well.

\subsection{Proposition} 
\label{pro for G}

\emph{
Let $C$ be a pre-operad. Then for all $h\t f\t g\in C^h\t C^f\t C^g$, 
the following composition relations hold:
}
$$
(h\o_i f)\o_j g=(-1)^{|f||g|}(h\o_{j-|f|}g)\o_i f, 
                       \qquad\text{if $i+f\leq j\leq|h|+|f|$}. 
$$

\subsection{Scope of a pre-operad}
\label{scope}

The {\it scope} of $(h\o_i f)\o_j g$ is given by
$$
0\leq i\leq|h|,\qquad 0\leq j\leq |f|+|h|.
$$
It follows from the defining relations of a pre-operad that
the scope is a disjoint union of 
\allowdisplaybreaks
\begin{align*}
B &\=\{(i,j)\in \NN\x\NN\,|\, 1\leq i\leq|h|  \,;\, 0   \leq j\leq i-1\},\\ 
A &\=\{(i,j)\in \NN\x\NN\,|\, 0\leq i\leq|h|  \,;\, i   \leq j\leq i+|f|\},\\ 
G &\=\{(i,j)\in \NN\x\NN\,|\, 0\leq i\leq|h|-1\,;\, i+f \leq j\leq |f|+|h|\}. 
\end{align*}
\allowdisplaybreaks
Note that the triangles $B$ and $G$ are symmetrically situated 
with respect to the parallelogram $A$ in the scope $BAG$.
The (recommended and impressive) picture is left for a reader as 
an exercise.

\subsection{Recapitulation}

The defining relations of a pre-operad can be easily rewritten as follows: 
$$
(h\o_i f)\o_j g=
\begin{cases} 
   (-1)^{|f||g|}(h\o_j g)\o_{i+|g|}f, &\text{if $(i,j)\in B$}\\ 
   h\o_i(f\o_{j-i}g),                 &\text{if $(i,j)\in A$}\\
   (-1)^{|f||g|}(h\o_{j-|f|}g)\o_i f, &\text{if $(i,j)\in G$}, 
\end{cases}
\medskip
$$
where we have included Proposition \ref{pro for G} as well.
The \emph{first} ($B$) and \emph{third} ($G$) parts of the relations 
turn out to be equivalent.

\section{Cup}
\subsection{Definition \cite{GGS92Am,CGS93}}

In a pre-operad $C$, let $\mu\in C^2$. Define
$\u\=\u_\mu\:C^f\t C^g\to C^{f+g}$ by
$$
f\u g\=(-1)^f(\mu\o_0 f)\o_f g,\qquad|\smile|=1,\qquad f\t g\in C^f\t C^g.
$$
The pair $\Cup\,C\=\{C,\u\}$ is called a $\u$-algebra of $C$.

\subsection{Example}
\emph{
For the endomorphism pre-operad (section \ref{HG}) $\EE_A$, 
one has 
}
$$
f\u g=(-1)^{fg}\mu\o(f\t g), 
      \qquad \mu\t f\t g\in \EE_A^2\t \EE_A^f\t \EE_A^g.
$$

\subsection{Proposition}
\label{cuppro}

\emph{
In a pre-operad $C$, one has
}
$$
\mu\o_0 f=(-1)^f f\u\I,
\quad \mu\o_1 f=-\I\u f,
\quad f\u g=-(-1)^{|f|g}(\mu\o_1 g)\o_0f.
$$
\begin{proof}
We have
$$
(-1)^f f\u\I=(-1)^{f+f}(\mu\o_0 f)\o_f\I=\mu\o_0 f,
           \quad -\I\u f=(\mu\o_0\I)\o_1 f=\mu\o_1 f.
$$
Also, calculate
\allowdisplaybreaks
\begin{align*}
f\u g &=(-1)^f(\mu\o_0 f)\o_f g
       =(-1)^{|f||g|+f}(\mu\o_{f-|f|}g)\o_0 f\\
      &=(-1)^{|f||g|+|f|+1}(\mu\o_1 g)\o_0 f
       =-(-1)^{|f|g}(\mu\o_1 g)\o_0 f,
\end{align*}
\allowdisplaybreaks
which was required as well.
\end{proof}

\subsection{Lemma}
\label{cup}
\emph{
In a pre-operad $C$, the following composition relations hold:
}
\smallskip
$$
(f\u g)\o_j h=
   \begin{cases}
     (-1)^{g|h|}(f\o_j h)\u g, &\text{if $0\leq j\leq|f|$}\\ 
     f\u(g\o_{j-f}h),          &\text{if $f\leq j\leq |g|+f$}.
\end{cases}
$$
\begin{proof}
Calculate, by using the defining relations of a pre-operad:
\allowdisplaybreaks
\begin{align*}
(f\u g)\o_j h &=(-1)^f[(\mu\o_0 f)\o_f g]\o_j h\\
              &= \begin{cases} 
(-1)^{f+|g||h|}[(\mu\o_0 f)\o_j h]\o_{f+|h|}g, &\text{if $0\leq j\leq|f|$}\\ 
(-1)^f(\mu\o_0 f)\o_f (g\o_{j-f}h),            &\text{if $f\leq j\leq |g|+f$}\\
                 \end{cases}\\
              &=\begin{cases}
(-1)^{f+|g||h|}[\mu\o_0(f\o_j h)]\o_{f+|h|}g, &\text{if $0\leq j\leq|f|$}\\ 
f\u(g\o_{j-f}h),                              &\text{if $f\leq j\leq |g|+f$}\\
                \end{cases}\\
              &=\begin{cases}
(-1)^{|f|+|h|+1+f+|g||h|}(f\o_j h)\u g, &\text{if $0\leq j\leq|f|$}\\ 
f\u(g\o_{j-f}h),                        &\text{if $f\leq j\leq |g|+f$}\\
                \end{cases}\\
              &=\begin{cases}
(-1)^{g|h|}(f\o_j h)\u g,               &\text{if $0\leq j\leq|f|$}\\ 
f\u(g\o_{j-f}h),                        &\text{if $f\leq j\leq |g|+f$},
                \end{cases}
\end{align*}
\allowdisplaybreaks
\par\medskip\noindent
which is the required formula.
\end{proof}

\section{Total composition and the Gerstenhaber identity} 
\subsection{Definition  \cite{GGS92Am,CGS93}}

In a pre-operad $C$, the \emph{total composition} 
$\bul\:C^f\t C^g\to C^{f+g-1}$ is defined by
$$ 
f\bul g\=\sum_{i=0}^{|f|}f\o_i g, 
     \qquad |\bul|=0, \qquad f\t g\in C^f\t C^g.
$$
The pair $\Com\,C\=\{C,\bul\}$ is called a \emph{composition algebra} of $C$.

\subsection{Theorem}
\label{right der}
\emph{
In a pre-operad $C$, one has
}
$$
(f\u g)\bul h=f\u(g\bul h)+(-1)^{|h|g}(f\bul h)\u g.
$$
\begin{proof}
Use Lemma \ref{cup}. Note that $|f\u g|=f+g-1$ and calculate,
\allowdisplaybreaks
\begin{align*}
(f\u g)\bul h &= \sum_{i=0}^{f+g-1}(f\u g)\o_i h 
 =\sum_{i=0}^{f-1}(f\u g)\o_i h+\sum_{i=f}^{f+g-1}(f\u g)\o_i h\\ 
&=(-1)^{|h|g}\sum_{i=0}^{|f|}(f\o_i h)\u g
  +\sum_{i=f}^{f+g-1}f\u(g\o_{i-f}h)\\ 
&=(-1)^{|h|g}(f\bul h)\u g+\sum_{i'=0}^{|g|}f\u(g\o_{i'}h)\\
&=(-1)^{|h|g}(f\bul h)\u g+f\u(g\bul h), 
\end{align*}
\allowdisplaybreaks
\noindent
which is the required formula.
\end{proof}

\subsection*{Remark} 

This theorem tells us that \emph{right} translations in $\Com\,C$ are 
(right) derivations of the $\u$-algebra. 
It may be anticipated from Theorem~5 of \cite{Ger} that the \emph{left}
translations in $\Com\,C$ are not derivations of the $\u$-algebra
(see section \ref{dev over ternary}).

\subsection{Associator}

Now, recall section \ref{scope} and note that
$$
(h\bul f)\bul g = \sum_{i=0}^{|h|}\sum_{j=0}^{|f|+|h|}(h\o_i f)\o_j g
                = (\sum_{(i,j)\in B}+\sum_{(i,j)\in A}+\sum_{(i,j)\in G})
                  (h\o_i f)\o_j g.
$$       
We can rearrange this double sum as follows. First note that
$$
\sum_{(i,j)\in A}(h\o_i f)\o_j g 
           =\sum_{i=0}^{|h|}\sum_{j=i}^{i+|f|}h\o_i(f\o_{j-i}g) 
           =\sum_{i=0}^{|h|}\sum_{j'=0}^{|f|}h\o_i(f\o_{j'}g) 
           =h\bul(f\bul g).
$$
So, an \emph{associator} is at hand,
$$
(h,f,g) \= (h\bul f)\bul g-h\bul(f\bul g)
         =  (\sum_{(i,j)\in B}+\sum_{(i,j)\in G})(h\o_i f)\o_j g.
$$

\subsection{Gerstenhaber braces \cite{Ger,CGS93}}

The Gerstenhaber ternary \emph{braces} $\{\cdot,\cdot,\cdot\}$ are 
defined as a double sum over the triangle $G$ by
$$
\{h,f,g\}\=\sum_{(i,j)\in G}(h\o_i f)\o_j g,
    \qquad |\{\cdot,\cdot,\cdot\}|=0,  \qquad h\t f\t g\in C^h\t C^f\t C^g.
$$

\subsection{Getzler identity \cite{Getzler}}
\label{Getzler}
\emph{
In a pre-operad $C$, one has
}
$$
(h,f,g)=\{h,f,g\}+(-1)^{|f||g|}\{h,g,f\}, 
            \qquad h\t f\t g\in C^h\t C^f\t C^g.
$$
\begin{proof}
First note that 
$$
\{h,f,g\}=\sum_{i=0}^{|h|-1}\sum_{j=i+f}^{|f|+|h|}(h\o_i f)\o_j g
         =\sum_{j=f}^{|f|+|h|}\sum_{i=0}^{j-f}(h\o_i f)\o_j g.
$$
By transposing the arguments, we have
$$
\{h,g,f\}=\sum_{j=g}^{|g|+|h|}\sum_{i=0}^{j-g}(h\o_i g)\o_j f
         =(-1)^{|f||g|}\sum_{j=g}^{|g|+|h|}\sum_{i=0}^{j-g}
                       (h\o_{j-|g|}f)\o_i g.
$$
Now introduce the new summation indices $i'$ and $j'$ by
$$
1\leq i'\=j-|g|\leq|h|,\qquad 0\leq j'\=i\leq i'-1.
$$
Then we have
$$
\{h,g,f\} =(-1)^{|f||g|}\sum_{i'=1}^{|h|}\sum_{j'=0}^{i'-1}
                      (h\o_{i'}f)\o_{j'}g
          =(-1)^{|f||g|}\sum_{(i,j)\in B}(h\o_i f)\o_j g,
$$
which proves the required formula.
\end{proof}

\subsection{Gerstenhaber identity \cite{Ger}}
\label{Gerst}
\emph{
In a pre-operad $C$, one has
}
$$
(h,f,g)=(-1)^{|f||g|}(h,g,f),\qquad h\t f\t g\in C^h\t C^f\t C^g.
$$
\begin{proof}
Use the Getzler identity \ref{Getzler}.
\end{proof}

\subsection{Remark}

Among others, we should like to call attention  particularly
to \ref{Gerst} and call it the \emph{Gerstenhaber identity}. 
This identity, first (form of) found in \cite{Ger}, is responsible for 
the Jacobi identity in $\Com^-C$ (defined below).
In Gerstenhaber's original terms from \cite{Ger}, it should be
called \emph{pre-Jacobi}.
We had quite the same motivation for \emph{Getzler's identity} \ref{Getzler},
thus it might be called a \emph{pre-Gerstenhaber} 
or \emph{pre-pre-Jacobi} identity. 

\subsection{Cup and braces}
\label{CaB}
\emph{
In a pre-operad $C$, one has
}
$$
f\u g=(-1)^f\{\mu,f,g\}, \qquad \mu\t f\t g\in C^2\t C^f\t C^g.
$$
\begin{proof}
Evidently, 
$$
\{\mu,f,g\}=\sum_{i=0}^{|\mu|-1}\sum_{j=i+f}^{|f|+|\mu|}(\mu\o_i f)\o_j g
           =\sum_{j=f}^{|f|+1}(\mu\o_0 f)\o_j g
           =(\mu\o_0 f)\o_f g
           =(-1)^f f\u g,
$$
which is the required formula.
\end{proof}

\section{Total composition and a graded Lie algebra}
\subsection{Commutator and Jacobian}

The \emph{commutator} $[\cdot,\cdot]$ and \emph{Jacobian} $J$ are defined 
in $\Com\,C$ in the conventional way by
\begin{align*}
[f,g]       &\=f\bul g-(-1)^{|f||g|}g\bul f=-(-1)^{|f||g|}[g,f],\\
J(f\t g\t h)&\=(-1)^{|f||h|}[[f,g],h]+(-1)^{|g||f|}[[g,h],f]
                                     +(-1)^{|h||g|}[[h,f],g].
\end{align*}
We denote the corresponding \emph{commutator algebra} of $C$ as 
$\Com^{-}C\=\{C,[\cdot,\cdot]\}$.

\subsection*{Remark}

The sign $^{-}$ stems from the definition of the commutator 
$[\cdot,\cdot]$. Quite a convenient notation
is $[f,g]_{\pm}\=f\bul g\pm(-1)^{|f||g|}g\bul f$.
We do not need $[\cdot,\cdot]_+$ in this paper, and so 
$[\cdot,\cdot]_-\=[\cdot,\cdot]$.

\subsection{Theorem \cite{Ger,GGS92Am,CGS93}} 

$\Com^-C$ \emph{is a graded Lie algebra.}
\begin{proof}
Indeed, in $\Com^-C$ the following identity holds:
\begin{align*}
J&(f\t g\t h)= (-1)^{|f||h|}[(f,g,h)-(-1)^{|g||h|}(f,h,g)]\\
             &+(-1)^{|g||f|}[(g,h,f)-(-1)^{|h||f|}(g,f,h)]
              +(-1)^{|h||g|}[(h,f,g)-(-1)^{|f||g|}(h,g,f)],
\end{align*}
so the Gerstenhaber identity implies the graded Jacobi identity: $J=0$.
\end{proof}

\noindent
In addition \cite{NiRi66}, one can easily check that
\begin{enumerate}
\item $[f,f]=0$,     if $|f|$ is \emph{even},
\item $[[f,f],f]=0$, if $|f|$ is \emph{odd}.
\end{enumerate}
The first item is evident from the definition of $[\cdot,\cdot]$. 
For the convenience of reader, let us check the second one. 
If $|f|$ is \emph{odd}, then
\begin{align*}
[[f,f],f]&=[f\bul f-(-1)^{|f||f|}f\bul f,f]=2[f\bul f,f]\\
         &=2[(f\bul f)\bul f-(-1)^{|f\bul f||f|}f\bul(f\bul f)]\\
         &=2(f,f,f)=2\{f,f,f\}+2(-1)^{|f||f|}\{f,f,f\}=0,
\end{align*}
where the Getzler identity was used inside the last row.

\subsection*{Remark}

The Jacobi identity implies that $3[[f,f],f]]=0$ for \emph{odd} $|f|$ and 
hence the restriction $\char\,K\neq3$ is imposed \cite{NiRi66} as a rule.
But via the Getzler identity one can avoid this unpleasant restriction.

\subsection{Remark}

This important theorem has been discussed several times.
Stasheff \cite{Sta92,Sta93} has proved this theorem for endomorphism 
pre-operads in \emph{homological} terms, by using \emph{coderivations} 
of a \emph{tensor coalgebra}. 
Getzler \cite{Getzler}, Getzler and Jones \cite{GetzJon94} have proved 
this theorem for graded endomorphism pre-operads.
Markl and Shnider \cite{MarklShn94,MarklShn96} have
proved this theorem for abstract graded pre-operads. 

\section{Cup and a pre-coboundary operator}
\subsection{Definition}

In a pre-operad $C$, define a \emph{pre-coboundary} operator $\de_\mu$ by 
$$
-\de_\mu f\=[f,\mu]\=\ad_\mu^{Right}f= f\bul\mu-(-1)^{|f|}\mu\bul f,
             \qquad \mu\t f\in C^2\t C^f.
$$

\subsection{Example}

In the Gerstenhaber theory \cite{Ger}, $C$ is an \emph{endomorphism}
pre-operad and $\de_\mu$ is the Hochschild \emph{coboundary operator} with 
the property $\de_\mu^2=0$, the latter is due to the \emph{associativity}
$\mu\bul\mu=0$ (see also section \ref{bul-sq}).

\subsection{Proposition}
\label{pre-coboundary}
\emph{
In a pre-operad $C$, one has
}
$$
-\de_\mu f=f\u\I+f\bul\mu+(-1)^{|f|}\,\I\u f,
\qquad \mu\t f\in C^2\t C^f.
$$

\subsection{Proposition}
\label{bul-sq}
\emph{
In a pre-operad $C$, one has $\de_\mu^2=-\de_{\mu\bul\mu}$. 
}
\begin{proof}
For $f\in C^f$ calculate 
\allowdisplaybreaks
\begin{align*}
\de^2_\mu f &= \de_\mu\de_\mu f=[[f,\mu],\mu]
             = [f\bul\mu-(-1)^{|f|}\mu\bul f,\mu]\\
            &=(f\bul\mu)\bul\mu-(-1)^{|f|}(\mu\bul f)\bul\mu
                               -(-1)^{|f|+1}\mu\bul(f\bul\mu)
                               -\mu\bul(\mu\bul f)\\
            &=f\bul(\mu\bul\mu)+(f,\mu,\mu)-(-1)^{|f|}(\mu,f,\mu)
                     -(\mu\bul\mu)\bul f+(\mu,\mu,f).
\end{align*}
\allowdisplaybreaks
Now note that Getzler's identity \ref{Getzler} implies $(f,\mu,\mu)=0$  
and the Gerstenhaber identity \ref{Gerst} implies 
$(\mu,\mu,f)=(-1)^{|f|}(\mu,f,\mu)$.
So, superfluous terms cancel out and we obtain 
$$
\de_\mu^2 f =f\bul(\mu\bul\mu)-(\mu\bul\mu)\bul f
            =[f,\mu\bul\mu]=\ad_{\mu\bul\mu}^{Right}f 
            =-\de_{\mu\bul\mu}f,
$$
which proves the required formula.
\end{proof}

\subsection*{Remark}

The standard proof of this theorem goes via the Jacobi identity. One has 
$$
2\de^2_\mu=[\de_\mu,\de_\mu]=[\ad_\mu,\ad_\mu]
          =\ad_{[\mu,\mu]}=2\ad_{\mu\bul\mu}=-2\de_{\mu\bul\mu},
$$
which means that $\char\,K\neq2$ is needed as a rule. 
By using the Getzler and Gerstenhaber identities one can avoid this 
restriction as well.

\section{Derivation deviation of $\de_\mu$ over total composition}
\subsection{Definition}

The \emph{derivation deviation} of $\de_\mu$ over $\bul$ is defined by 
$$
\dev_\bul\de_\mu(f\t g)
   \=\de_\mu(f\bul g)-f\bul\de_\mu g-(-1)^{|g|}\de_\mu f\bul g.
$$

\subsection{Theorem \cite{Ger}}
\label{second main}
\emph{
In a pre-operad $C$, one has
}
$$
(-1)^{|g|}\dev_\bul\de_\mu(f\t g)=f\u g-(-1)^{fg}g\u f,
\qquad \mu\t f\t g\in C^2\t C^f\t C^g.
$$
\begin{proof}
First use the definitions of $\de_\mu$ and $[\cdot,\cdot]$:
\allowdisplaybreaks
\begin{align*}
\dev_\bul\de_\mu(f\t g)
     \=&\,\,\de_\mu(f\bul g)-f\bul\de_\mu g-\,(-1)^{|g|}\de_\mu f\bul g\\
      =&-[f\bul g,\mu]+f\bul[g,\mu]+(-1)^{|g|}[f,\mu]\bul g\\
      =&-(f\bul g)\bul\mu+(-1)^{|f|+|g|}\mu\bul(f\bul g)+f\bul(g\bul\mu)\\
       &-(-1)^{|g|}f\bul(\mu\bul g)+(-1)^{|g|}(f\bul\mu)\bul g 
        -(-1)^{|g|+|f|}(\mu\bul f)\bul g\\  
      =&-(f,g,\mu)-(-1)^{|f|+|g|}(\mu,f,g)+(-1)^{|g|}(f,\mu,g).
\end{align*}
\allowdisplaybreaks
Now, note that $(f,g,\mu)=(-1)^{|g|}(f,\mu,g)$
and use Getzler's identity \ref{Getzler} with Proposition \ref{CaB}. 
So it follows that
\allowdisplaybreaks
\begin{align*}
\dev_\bul\de_\mu(f\t g)
    &=-(-1)^{|f|+|g|}(\mu,f,g)\\
    &=-(-1)^{|f|+|g|}\{\mu,f,g\}-(-1)^{|f|+|g|+|f||g|}\{\mu,g,f\}\\
    &=-(-1)^{|f|+|g|+f}f\u g-(-1)^{|f|+|g|+|f||g|+g}g\u f\\
    &=\,\,(-1)^{|g|}[f\u g-(-1)^{fg}g\u f],
\end{align*}
\allowdisplaybreaks
which is the required formula.
\end{proof}

\subsection{Remark}

An alternative proof (Gertenhaber's method) of Theorem \ref{second main} 
is presented in the next section, read also sections \ref{aux} and 
\ref{rem aux} for motivations. 

\section{Revival of the Gerstenhaber method 
    (second proof of Theorem \ref{second main})}
\label{sec main}

\noindent
In this, quite a didactic section, we should like to illustrate in detail
the essence of Gerstenhaber's method \cite{Ger} when calculating 
\emph{derivation deviations} of the \emph{pre-coboundary} operator.

\subsection{Auxiliary variables}
\label{aux var}

In a pre-operad $C$, for $f\t g\in C^f\t C^g$ define (cf.~\cite{Ger}) 

\allowdisplaybreaks
\begin{align*}
\la_{i+1}&\=-(-1)^{|f|+|g|}\,\I\u(f\o_i g)
            -(-1)^{|g|}\sum_{j=0}^{i-1}(f\o_j\mu)\o_{i+1}g
            +(-1)^{|g|}f\o_i(\I\u g), \\
\la'_{i+1}&\= f\o_i(g\u\I)
             -(-1)^{|g|}\sum_{j=i+1}^{|f|}(f\o_j\mu)\o_ig
             -(f\o_i g)\u\I,  \qquad\quad\,\, 0\leq i\leq|f|.
\end{align*}
\allowdisplaybreaks

\subsection{Lemma}
\label{first}
\emph{
In a pre-operad $C$, one has
}
$$
\de_\mu(f\o_i g)-f\o_i\de_\mu g=\la_{i+1}+\la'_{i+1},\qquad 0\leq i\leq|f|. 
$$ 
\begin{proof} 
See Appendix A. 
\end{proof}

\subsection{Lemma}
\label{second}
\emph
{In a pre-operad $C$, one has
}
$$
(-1)^{|g|}\de_\mu f\o_i g=\la_i+\la'_{i+1},\qquad 0\leq i\leq f,
$$ 
\emph{by definition for $\la_0$ and $\la'_{f+1}$.}
\begin{proof}
See Appendix B. 
\end{proof}

\subsection{Gerstenhaber's method 
            {\rm(via a proof of Theorem \ref{second main})}}

By using Lemma \ref{first} and Lemma \ref{second} we have
\allowdisplaybreaks
\begin{align*}
\dev_\bul\de_\mu(f\t g)
  &= \sum_{i=0}^{|f|}[\de_\mu(f\o_i g)-f\o_i\de_\mu g]
         -(-1)^{|g|}\sum_{i=0}^{f}\de_\mu f\o_i g\\
  &= \sum_{i=0}^{|f|}(\la_{i+1}+\la'_{i+1})
                       -\sum_{i=0}^{f}(\la_i+\la'_{i+1})\\
  &= \sum_{i=0}^{|f|}\la_{i+1}-\sum_{i=0}^{f}\la_i
    +\sum_{i=0}^{|f|}\la'_{i+1}-\sum_{i=0}^{f}\la'_{i+1}
  =\framebox{$-\la_0-\la'_{f+1}$}\\
  &=\la_f+\la'_1-(-1)^{|g|}(\de_\mu f\o_0 g+\de_\mu f\o_f g).
\end{align*}
\allowdisplaybreaks
Calculate
\allowdisplaybreaks
\begin{align*}
\la_f+\la'_1=&-(-1)^{|f|+|g|}\,\I\u(f\o_{|f|}g) 
             -(-1)^{|g|}\sum_{j=0}^{|f|-1}(f\o_j\mu)\o_f g\\
             &-(-1)^{|g|}(f\o_{|f|}\mu)\o_f g
              +(-1)^{|g|}(f\o_{|f|}\mu)\o_f g\\
             &+(-1)^{|g|}f\o_{|f|}(\I\u g) 
             + f\o_0(g\u I)\\
            &-(-1)^{|g|}\sum_{j=1}^{|f|}(f\o_j\mu)\o_0 g
             -(-1)^{|g|}(f\o_0\mu)\o_0 g\\
            &+(-1)^{|g|}(f\o_0\mu)\o_0 g
             -(f\o_0 g)\u\I\\
\medskip
           =&-(-1)^{|f|+|g|}\,\I\u(f\o_{|f|}g)
             -(-1)^{|g|} (f\bul\mu)\o_f g
             +(-1)^{|g|}(f\o_{|f|}\mu)\o_f g\\
            &+(-1)^{|g|}f\o_{|f|}(\I\u g)
             +f\o_0(g\u I)
             -(-1)^{|g|}(f\bul\mu)\o_0 g\\
            &+(-1)^{|g|}(f\o_0\mu)\o_0 g
             -(f\o_0 g)\u\I\\
\medskip
           =&-(-1)^{|f|+|g|}\,\I\u(f\o_{|f|}g)
             +(-1)^{|f|+|g|}(\I\u f)\o_f g
             +(-1)^{|g|}\de f\o_f g\\
            &+(-1)^{|g|}(f\u\I)\o_f g
             +(-1)^{|g|}(f\o_{|f|}\mu)\o_f g
             +(-1)^{|g|}f\o_{|f|}(\I\u g)\\
            &+f\o_0(g\u I)
             +(-1)^{|f|+|g|}(\I\u f)\o_0 g
             +(-1)^{|g|}\de f\o_0 g\\
            &+(-1)^{|g|}(f\u\I)\o_0 g
             +(-1)^{|g|}(f\o_0\mu)\o_0 g
             -(f\o_0 g)\u\I.
\end{align*}
\allowdisplaybreaks
Now, use composition relations to note that
\allowdisplaybreaks
\begin{gather*}
(-1)^{|f|+|g|}(\I\u f)\o_f g   = (-1)^{|f|+|g|}\,\I\u(f\o_{|f|}g),\\
(-1)^{|g|} (f\u\I)\o_f g       = (-1)^{|g|}f\u(\I\o_0 g)
                               = (-1)^{|g|}f\u g,\\
(-1)^{|g|}(f\o_{|f|}\mu)\o_f g = (-1)^{|g|}f\o_{|f|}(\mu\o_1 g)
                               = -(-1)^{|g|}f\o_{|f|}(\I\u g),\\
(-1)^{|f|+|g|}(\I\u f)\o_0 g   = (-1)^{|f|+|g|+f|g|}(\I\o_0 g)\u f
                               = (-1)^{g|f|}g\u f,\\
(-1)^{|g|}(f\u\I)\o_0 g        = (-1)^{|g|+|g|}(f\o_0 g)\u\I
                               = (f\o_0 g)\u\I,\\
(-1)^{|g|}(f\o_0\mu)\o_0 g     = (-1)^{|g|}f\o_0(\mu\o_0 g)
                               =-f\o_0(g\u\I).
\end{gather*}
\allowdisplaybreaks
So, many terms cancel out and we have
\begin{align*}
\la_f+\la'_1&=(-1)^{|g|}f\u g-(-1)^{g|f|-1}g\u f
                          +(-1)^{|g|}(\de f\o_f g+\de f\o_0 g)\\ 
        &=(-1)^{|g|}[f\u g-(-1)^{fg}g\u f]
                          +(-1)^{|g|}(\de f\o_f g+\de f\o_0 g). 
\end{align*}
By substituting this into $\dev_\bul\de_\mu(f\t g)$, we obtain the required formula.
\qed

\subsection{Remark}
\label{rem aux}

At this point, it may be asked that why does one need such a lengthy proof,
if a much shorter one is at hand from section \ref{second main}.
An idea of the above proof can be stated in a few words as follows.

\subsection*{\emph{Observation}}
\emph{
The r.~h.~s. of formula \ref{second main} is the sum of the 
boundary terms $\la_0$ and $\la'_{f+1}$.
}

\smallskip
Lemmas \ref{first} and \ref{second} with the auxiliary variables 
$\la_i$ and $\la'_i$ serve only as a remedy for a time.
The essential point, discovered by Gerstenhaber \cite{Ger}, is that
the above \emph{auxiliary variables} method can be adapted (generalized) 
also for a considerably more complicated situation described concisely 
in the next section, where any shorter proof  
(like the proof presented in the section \ref{second main})
is not known.

\section{Discussion: Derivation deviation over ternary braces}
\label{dev over ternary}

\subsection{Definition}
The \emph{derivation deviation} of $\de\=\de_\mu$ over 
$\{\cdot,\cdot,\cdot\}$ is defined by 
$$
\dev_{\{\cdot,\cdot,\cdot\}}\de\,(h\t f\t g)
    \=\de\{h,f,g\}-\{h,f,\de g\}-(-1)^{|g|}\{h,\de f,g\}
      -(-1)^{|g|+|f|}\{\de h,f,g\}.
$$

\subsection{Theorem {\rm (cf.~\cite{Ger,VorGer})}}
\label{third main}
\emph{
In a pre-operad $C$, one has
}
$$
(-1)^{|g|}\dev_{\{\cdot,\cdot,\cdot\}}\de\,(h\t f\t g)=
    (h\bul f)\u g+(-1)^{|h|f}f\u(h\bul g)-h\bul(f\u g).
$$

\subsection*{Remark instead of proof}
\label{aux}

The detailed proof of this theorem will be exposed in \cite{KP}. 
As in the above example of the method,
the proof can be performed via \emph{auxiliary variables}.
Originally, Gerstenhaber's proof of the Theorem 5 in \cite{Ger} 
was presented in the \emph{endomorphism} pre-operad terms.
It was announced in \cite{Ger68,GGS92Am,CGS93} that the Theorem
also holds in an abstract setting, i.~e. for abstract pre-operads.
However, the adaptation of Gerstenhaber's original proof for an
\emph{abstract} pre-operad is not trivial.
In short, an idea of the proof is as follows.

\smallskip
\noindent
{\bf\emph{Observation in advance.}}
\emph{
The r.~h.~s. of formula \ref{third main} turns out to be a sum of the 
boundary terms of the (suitably defined) auxiliary variables. 
The sum is over the boundary of a truncated triangle
enveloping the Gerstenhaber triangle $G$.
}

\smallskip
In this sense, the \emph{simplicial} structure of the defining $\o_i$ 
relations turns up. 
Recall, that already in the early days of comp calculus,
$\u$ and $\bul$ were recognized as Steenrod operations \cite{Steenrod}.

\subsection{Theorem {\rm(cf.~\cite{Ger,VorGer})}}
\emph{
In a pre-operad $C$, one has
}
$$
(-1)^{|g|}\dev_{\{\cdot,\cdot,\cdot\}}\de\,(h\t f\t g)=
    [h,f]\u g+(-1)^{|h|f}f\u[h,g]-[h,f\u g].
$$
\begin{proof} Combine the previous Theorem \ref{third main} with Theorem
\ref{right der}.
\end{proof}

\subsection*{Remark}
A well known first form of this theorem, found by Gerstenhaber \cite{Ger} 
for the Hochschild complex, can be seen as a starting point of the modern
\emph{mechanical mathematics} based nowadays on the pioneering concept of 
(\emph{homotopy} \cite{GerVor94,VorGer}) \emph{Gerstenhaber algebra}. 

\begin{center}
\section*{Appendix A}
\end{center}

\subsection*{Proof of Lemma \ref{first}}
First note that
$$
\de_\mu(f\o_i g) =-(-1)^{|f|+|g|}\,\I\u(f\o_i g)
                  -\sum_{j=0}^{|f|+|g|}(f\o_i g)\o_j\mu
                  -(f\o_i g)\u\I.
$$
By substituting here
$$
(f\o_i g)\o_j\mu=
\begin{cases}
     (-1)^{|g|}(f\o_j\mu)\o_{i+1}g,   &\text{if $0\leq j\leq i-1$}      \\
     f\o_i(g\o_{j-i}\mu),             &\text{if $i\leq j\leq i+|g|$}    \\
     (-1)^{|g|}(f\o_{j-|g|}\mu)\o_ig, &\text{if $i+g\leq j\leq |g|+|f|$}\\
\end{cases}
$$
we have
\allowdisplaybreaks
\begin{align*}
\de_\mu(f\o_i g)= &-(-1)^{|f|+|g|}\,\I\u(f\o_i g)
                   -(-1)^{|g|}\sum_{j=0}^{i-1}(f\o_j\mu)\o_{i+1}g\\
                  &-\sum_{j=i}^{i+|g|}f\o_i(g\o_{j-i}\mu)
                   -(-1)^{|g|}\sum_{j=i+g}^{|g|+|f|}(f\o_{j-|g|}\mu)\o_i g
                   -(f\o_i g)\u\I.
\end{align*}
\allowdisplaybreaks
Now, note that
$$
\sum_{j=i}^{i+|g|}f\o_i(g\o_{j-i}\mu)=\sum_{j=0}^{|g|}f\o_i(g\o_j\mu)
$$
and use 
$$
\sum_{j=0}^{|g|}g\o_j\mu=g\bul\mu
                        =-(-1)^{|g|}\,\I\u g-\de_\mu g-g\u\I.
$$
Then we obtain
\allowdisplaybreaks
\begin{align*}
\de_\mu(f\o_i g)=&-(-1)^{|f|+|g|}\,\I\u(f\o_i g)
                   -(-1)^{|g|}\sum_{j=0}^{i-1}(f\o_j\mu)\o_{i+1}g\\
                  &+(-1)^{|g|}f\o_i(\I\u g)
                   +f\o_i\de_\mu g
                   +f\o_i(g\u\I)\\
                  &-(-1)^{|g|}\sum_{j=i+1}^{|f|}(f\o_j\mu)\o_i g
                   -(f\o_i g)\u\I \\ 
                 =&+\la_{i+1}+f\o_i\de_\mu g+\la'_{i+1},
\end{align*}
\allowdisplaybreaks
which proves the required formula.
\qed

\begin{center}
\section*{Appendix B} 
\end{center}

\subsection*{Proof of Lemma \ref{second}}
First note that 
\allowdisplaybreaks
\begin{align*}
\la_i+\la'_{i+1}=&-(-1)^{|f|+|g|}\,\I\u(f\o_{i-1}g)
                  -(-1)^{|g|}\sum_{j=0}^{i-2}(f\o_j\mu)\o_i g\\
                 &+(-1)^{|g|}f\o_{i-1}(\I\u g)
                  +f\o_i(g\u\I)\\
                 &-(-1)^{|g|}\sum_{j=i+1}^{|f|}(f\o_j\mu)\o_i g
                  -(f\o_i g)\u\I.
\end{align*}
\allowdisplaybreaks
We must compare it term by term with 
\allowdisplaybreaks
\begin{align*}
\de_\mu f\o_i g=&-[(-1)^{|f|}\,\I\u f
                 -\sum_{j=0}^{|f|}(f\o_j\mu)
                 -f\u\I]\o_i g\\
               =&-(-1)^{|f|}(\I\u f)\o_i g
                 -\sum_{j=0}^{|f|}(f\o_j\mu)\o_i g
                 -(f\u\I)\o_i g\\
               =&-(-1)^{|f|}(\I\u f)\o_i g
                 -\sum_{j=0}^{i-2}(f\o_j\mu)\o_i g
                 -(f\o_{i-1}\mu)\o_i g\\
                &-(f\o_i\mu)\o_i g
                 -\sum_{j=i+1}^{|f|}(f\o_j\mu)\o_i g
                 -(f\u\I)\o_i g.
\end{align*}
\allowdisplaybreaks
Now, use composition relations to note that
\allowdisplaybreaks
\begin{gather*}
(\I\u f)\o_i g=\I\u(f\o_{i-1}g),\\
(f\o_{i-1}\mu)\o_i g=f\o_{i-1}(\mu\o_1 g)=-f\o_{i-1}(\I\u g),\\
(f\o_i\mu)\o_i g=f\o_i(\mu\o_0 g)=(-1)^g f\o_i(g\u\I),\\
(f\u\I)\o_i g=(-1)^{|g|}(f\o_i g)\u\I,
\end{gather*}
\allowdisplaybreaks
which proves the required formula. \qed

\medskip
\section*{\centerline{ACKNOWLEDGEMENTS}}
\noindent
We thank Martin Markl and the referee for thorough reading of the original 
version, helpful remarks and suggestions.

\end{document}